\documentclass[11pt,letterpaper]{article}
\usepackage[T1]{fontenc}
\usepackage{lmodern}
\usepackage{amsmath,amsthm,amssymb}
\usepackage[margin=1.05in]{geometry}
\usepackage{microtype}
\usepackage{xcolor}
\usepackage{enumitem}
\usepackage{tikz}
\usetikzlibrary{arrows.meta,positioning,calc,backgrounds}
\usepackage{fancyhdr}
\usepackage{hyperref}
\definecolor{paperlink}{RGB}{0,0,114}
\setlist[enumerate]{leftmargin=2.5em,itemsep=.3em,topsep=.4em}
\newtheorem{theorem}{Theorem}[section]
\newtheorem{lemma}[theorem]{Lemma}
\newtheorem{proposition}[theorem]{Proposition}
\newtheorem{corollary}[theorem]{Corollary}
\newtheorem{claim}{Claim}
\newtheorem{problem}{Problem}[section]
\numberwithin{equation}{section}
\DeclareMathOperator{\diam}{diam}
\newcommand{\strictto}{\mathrel{\mapsto}}
\newcommand{\nonadj}{\mathrel{\not\sim}}
\tikzset{vertex/.style={circle,draw,fill=white,inner sep=1.8pt},arc/.style={-{Stealth[length=2mm]},semithick},segment/.style={-{Stealth[length=2mm]},thick},crossarc/.style={-{Stealth[length=2mm]},thick,draw=paperlink}}

\definecolor{cycleC}{HTML}{202020}
\definecolor{cycleE}{HTML}{202020}
\definecolor{insertarc}{HTML}{202020}
\definecolor{deletearc}{HTML}{505050}
\tikzset{cv/.style={circle,draw=cycleC,fill=white,line width=.8pt,inner sep=2.4pt},ev/.style={circle,draw=cycleE,fill=white,line width=.8pt,inner sep=2.4pt},keepC/.style={-{Stealth[length=2.2mm]},draw=cycleC,line width=1.3pt},keepE/.style={-{Stealth[length=2.2mm]},draw=cycleE,line width=1.3pt},added/.style={-{Stealth[length=2.2mm]},draw=insertarc,line width=1.3pt,dash pattern=on 5pt off 2pt on 1pt off 2pt},removed/.style={-{Stealth[length=2.2mm]},draw=deletearc,line width=1pt,densely dashed},panel/.style={anchor=west,font=\small\bfseries,text=cycleC},badge/.style={circle,fill=insertarc,text=white,font=\scriptsize\bfseries,inner sep=2pt}}

\title{Hamiltonicity and cycle factors in $k$-quasi-transitive digraphs of large diameter}
\author{Hao Zhang\;\quad Panpan Cheng\thanks{Corresponding author: Panpan Cheng (ppancheng@126.com)}\\
	\small{School of Mathematics, Shandong University, Jinan 250100, China}
	}
\date{}
\date{}
\begin{document}
\maketitle
\begin{abstract}
A digraph is $k$-quasi-transitive if the endpoints of every directed path of length $k$ are adjacent. In this paper, we prove that every strong non-bipartite $k$-quasi-transitive digraph of diameter at least $k+2$ is Hamiltonian when $k\ge7$ is odd. Moreover, combined with the known even case, for $k\ge6$, we also prove that a strong $k$-quasi-transitive digraph of diameter at least $k+2$ is Hamiltonian if and only if it has a cycle factor. These results answer two problems posed by Alva-Samos and Hern\'andez-Cruz affirmatively in [Journal of Graph Theory \textbf{96} (2021) 310--325].
\end{abstract}
\noindent\textbf{Keywords.} $k$-quasi-transitive digraph; Hamiltonian cycle; cycle factor; directed separator.

\section{Introduction}
A digraph is \emph{$k$-quasi-transitive} if the endpoints of every directed path of length $k$ are adjacent. For $k=2$, this is the classical class of quasi-transitive digraphs. Bang-Jensen and Huang~\cite{BH} characterized the Hamiltonian members of this class, and Bang-Jensen and Gutin~\cite{BG} gave polynomial algorithms for finding vertex-heaviest paths and cycles with nonnegative vertex weights.

A \emph{cycle factor} is a collection of vertex-disjoint directed cycles covering all vertices. Every Hamiltonian digraph is strong and has a cycle factor, but the converse fails even for quasi-transitive digraphs, as examples in~\cite{BH} show. Bang-Jensen and Nielsen~\cite{BN} gave a polynomial procedure for reducing a cycle factor to an irreducible one, but showed that the resulting factor need not be minimum. A cycle-factor characterization of Hamiltonicity therefore requires a structural argument ensuring that the reduction can reach a single cycle.

For bipartite tournaments, H\"aggkvist and Manoussakis~\cite{HM} characterized Hamiltonicity in terms of strong connectivity and a cycle factor. Cano, Galeana-S\'anchez and Goldfeder~\cite{CGG} obtained a corresponding characterization for $\mathcal P$-compositions of directed cycles whose components are strong semicomplete bipartite digraphs. Such conclusions depend on the structure of the class: Bang-Jensen, Gutin and Huang~\cite{BGH} showed that no fixed degree of strong connectivity, together with a cycle factor, guarantees Hamiltonicity for multipartite tournaments.

Hern\'andez-Cruz and Galeana-S\'anchez~\cite{HG} investigated kernels and reachability in $k$-transitive and $k$-quasi-transitive digraphs. For $k=3$, Galeana-S\'anchez, Goldfeder and Urrutia~\cite{GGU} obtained a complete description of the strong digraphs, involving semicomplete digraphs, semicomplete bipartite digraphs and an additional infinite family. For larger $k$, structural results are available when the diameter is at least $k+2$. A shortest path of length $k+2$ then constrains both the arcs among its vertices and the neighborhoods of the remaining vertices.

For even $k\ge4$, Wang and Zhang~\cite{WZ} proved that every strong $k$-quasi-transitive digraph of diameter at least $k+2$ has a Hamiltonian path. Hern\'andez-Cruz and Ju\'arez-Camacho~\cite{HJ} subsequently proved pancyclicity for even $k>4$, settling Hamiltonicity in that range. For odd $k\ge5$, Wang and Zhang~\cite{WZodd} described the semicomplete or semicomplete bipartite structure induced by a shortest path of length $k+2$, together with the possible structures of its complement. 
Alva-Samos and Hern\'andez-Cruz~\cite{AH} established structural results under the same diameter hypothesis and posed the following problem.

%Alva-Samos and Hern\'andez-Cruz~\cite{AH} established structural results under the same diameter hypothesis and posed questions about Hamiltonian cycles. In particular, they ask whether the strong non-bipartite digraphs are Hamiltonian. Our first theorem answers this question for odd $k\ge7$.

%\begin{theorem}\label{thm:odd}
%Let $k\ge7$ be odd. Every strong non-bipartite $k$-quasi-transitive digraph $D$ with $\diam(D)\ge k+2$ is Hamiltonian.
%\end{theorem}

%Moreover, Alva-Samos and Hern\'andez-Cruz \cite{AH} also  ask whether a cycle factor suffices for Hamiltonicity. The next theorem answers this question affirmatively; its even case follows from~\cite{HJ}.

%\begin{theorem}\label{thm:factor}
%Let $k\ge6$, and let $D$ be a strong $k$-quasi-transitive digraph with $\diam(D)\ge k+2$. Then $D$ is Hamiltonian if and only if it has a cycle factor.
%\end{theorem}

\begin{problem}[{\cite{AH}}]\label{prob:odd}
Let $k\ge7$ be odd. Is every strong non-bipartite $k$-quasi-transitive digraph $D$ with $\diam(D)\ge k+2$ Hamiltonian?
\end{problem}

Our first theorem answers Problem~\ref{prob:odd} affirmatively.

\begin{theorem}\label{thm:odd}
Let $k\ge7$ be odd. Every strong non-bipartite $k$-quasi-transitive digraph $D$ with $\diam(D)\ge k+2$ is Hamiltonian.
\end{theorem}

They also posed the following problem concerning cycle factors, which we state with the necessary strong connectivity assumption made explicit.

\begin{problem}[{\cite{AH}}]\label{prob:factor}
Let $k\ge6$. Is a strong $k$-quasi-transitive digraph $D$ with $\diam(D)\ge k+2$ Hamiltonian if and only if it has a cycle factor?
\end{problem}

The following theorem answers Problem~\ref{prob:factor} affirmatively. The case of even $k$ follows from~\cite{HJ}.

\begin{theorem}\label{thm:factor}
Let $k\ge6$, and let $D$ be a strong $k$-quasi-transitive digraph with $\diam(D)\ge k+2$. Then $D$ is Hamiltonian if and only if it has a cycle factor.
\end{theorem}

%The difficulty in the odd bipartite case lies in the interaction between a shortest path of length $k+2$ and the vertices outside it. The results of~\cite{AH} give complete adjacency between opposite-part vertices within each set, while missing adjacencies may remain between the sets. We prove that a cycle containing an outside vertex and meeting the shortest path must also contain every opposite-part vertex of that path nonadjacent to the outside vertex. This permits all factor cycles meeting the path to be merged. We then absorb the remaining cycles that contain vertices with missing adjacencies to the path. Vertices on distinct remaining cycles are then adjacent whenever they lie in opposite parts, so a final merging argument produces a Hamiltonian cycle.

The paper is organized as follows. Section~\ref{sec:structure} presents some preliminary results. In Section~\ref{sec:non-bipartite}, we prove Theorem~\ref{thm:odd}. Sections~\ref{sec:separation} and~\ref{sec:merging} establish several lemmas needed to prove Theorem~\ref{thm:factor}. Finally, we prove Theorem~\ref{thm:factor} in Section~\ref{sec:factor}.

\section{Preliminaries}\label{sec:structure}
All digraphs are finite, without loops or parallel arcs; opposite arcs are allowed. We write $x\to y$ when $xy$ is an arc. Distinct vertices $x,y$ are \emph{adjacent} if at least one of $x\to y$ and $y\to x$ holds. We write $x\nonadj y$ otherwise. For disjoint vertex sets $S,T$, the notation $S\to T$ means that every vertex of $S$ dominates every vertex of $T$. We write $S\strictto T$ if, in addition, no arc is directed from $T$ to $S$. Singletons may be replaced by their vertices in this notation.

A digraph is \emph{semicomplete} if every two distinct vertices are adjacent. A digraph is \emph{bipartite} if its vertex set has a partition $(X,Y)$ with no arcs within either part. Such a digraph is \emph{semicomplete bipartite} if every pair in $X\times Y$ is adjacent. For induced subdigraphs, one of the inherited parts may be empty. The induced subdigraph on a vertex set $S$ is denoted by $D[S]$, and $D-v$ denotes the deletion of a vertex $v$.

All paths and cycles are directed and simple, and their lengths count arcs. Cycles of length two are permitted. If $C$ is a directed cycle and $x,y$ are distinct vertices on it, then $xCy$ denotes the directed segment from $x$ to $y$ along $C$; the same notation is used for a directed path when the vertices occur in the indicated order. For a vertex $x$ on a directed cycle, $x^+$ denotes its successor on that cycle. A Hamiltonian cycle visits every vertex.

The distance $d_D(u,v)$ (or $d(u,v)$ for short) is the minimum length of a $u$--$v$ path in $D$, or infinity if no such path exists. We put $\diam(D)=\max_{u,v\in V(D)}d(u,v)$. A digraph is \emph{strong} if every vertex can reach every other vertex. If $D$ is strong and $\diam(D)\ge k+2$, an initial segment of a sufficiently long shortest path gives a path
\begin{equation}\label{eq:geodesic}
P=v_0v_1\cdots v_{k+2}\text{ and } d(v_0,v_{k+2})=k+2\,.
\end{equation}
Fix such a path and put $A=V(P)$ and $B=V(D)\setminus A$. We call $A$ the \emph{core}. Every subpath of $P$ is shortest, so $d(v_i,v_j)=j-i$ whenever $i<j$. In particular, there is no arc $v_i\to v_j$ with $j\ge i+2$.
Moreover, let $A_p=\{v_i:0\le i\le k+2,\ i\equiv p\pmod2\}$ for $p\in\{0,1\}$. For a vertex $z$, let $N(z)$ denote the set of vertices adjacent to $z$. 

\subsection{Lemmas}
We collect the following lemmas of Alva-Samos and Hern\'andez-Cruz \cite{AH} to prove Propositions~\ref{prop:semi} and~\ref{prop:bip}. Throughout Lemmas~\ref{src:4}--\ref{src:19}, let $k\ge7$ be odd, let $D$ be a strong $k$-quasi-transitive digraph, and let $P=v_0v_1\cdots v_{k+2}$ be a directed path with $d(v_0,v_{k+2})=k+2$. Retain the notation $A$, $B$, $A_p$ and $N(z)$ defined above. These common assumptions and definitions apply to all ten lemmas below.

\begin{lemma}[{\cite{AH}}]\label{src:4}
For any indices $i,j\in\{0,1,\ldots,k+2\}$ with $i\not\equiv j\pmod2$, the vertices $v_i$ and $v_j$ of $P$ are adjacent.
\end{lemma}

\begin{lemma}[{\cite{AH}}]\label{src:7}
Every vertex $z\in B$ is adjacent to some $v_i$ with $0\le i\le k+2$.
\end{lemma}

\begin{lemma}[{\cite{AH}}]\label{src:8}
Let $z\in B$. If $z\to v_i$ with $2\le i\le k+2$, then $z$ is adjacent to $v_{i-2}$ and $z\to v_j$ for every $0\le j\le i-4$ with $j\equiv i\pmod2$. If $v_i\to z$ with $0\le i\le k$, then $z$ is adjacent to $v_{i+2}$ and $v_j\to z$ for every $i+4\le j\le k+2$ with $j\equiv i\pmod2$.
\end{lemma}

\begin{lemma}[{\cite{AH}}]\label{src:9}
For $z\in B$, adjacency to any one vertex of $\{v_0,v_2,v_{k-1},v_{k+1}\}$ is equivalent to adjacency to each of them. The same holds for $\{v_1,v_3,v_k,v_{k+2}\}$.
\end{lemma}

\begin{lemma}[{\cite{AH}}]\label{src:10}
Suppose that $z\in B$ has a neighbor in $A_p$, where $p\in\{0,1\}$. Then either $z\strictto A_p$, or $A_p\strictto z$, or $z$ is adjacent to every vertex of $A_p$ with both an outneighbor and an inneighbor there, or its outneighbors in $A_p$ are exactly the first $r$ vertices and its inneighbors in $A_p$ exactly the last $s$ vertices, in increasing index order, where $r,s\ge2$ and $r+s<|A_p|$.
\end{lemma}

\begin{lemma}[{\cite{AH}}]\label{src:13}
If $D[A]$ is semicomplete, each $z\in B$ satisfies one of the following: $z\strictto A$; $A\strictto z$; $z$ is adjacent to all of $A$ and has both an outneighbor and an inneighbor in $A$; or $\{v_s,\ldots,v_{k+2}\}\strictto z\strictto\{v_0,\ldots,v_r\}$ for some $3\le r<s-1\le k-2$.
\end{lemma}

\begin{lemma}[{\cite{AH}}]\label{src:14}
If there exist indices $0\le i<j\le k+2$ with $i\equiv j\pmod2$ such that $v_i$ and $v_j$ are adjacent, then $D[A]$ is semicomplete.
\end{lemma}

\begin{lemma}[{\cite{AH}}]\label{src:15}
If some $z\in B$ has a neighbor in each of $A_0$ and $A_1$, then $D[A]$ is semicomplete.
\end{lemma}

\begin{lemma}[{\cite{AH}}]\label{src:16}
Distinct $x,y\in B$ are adjacent whenever $N(x)\cap A_0\ne\varnothing$ and $N(y)\cap A_1\ne\varnothing$.
\end{lemma}

\begin{lemma}[{\cite{AH}}]\label{src:19}
If adjacent vertices $u,v\in B$ both have a neighbor in $A_p$ for some $p\in\{0,1\}$, then $D[A]$ is semicomplete.
\end{lemma}

\subsection{Propositions}

\begin{proposition}\label{prop:semi}
Suppose that $k\ge7$ is odd and $D$ is strong, non-bipartite, and $k$-quasi-transitive. Then $D[A]$ and $D[B]$ are semicomplete. If $z\in B$ is nonadjacent to some vertex of $A$, then $v_{k+2}\to z\to v_0$.
\end{proposition}
\begin{proof}
We first prove that $D[A]$ is semicomplete. Suppose otherwise. An adjacent pair of core vertices with indices of the same parity would force $D[A]$ to be semicomplete by~Lemma~\ref{src:14}. Thus $A_0$ and $A_1$ are independent. By Lemma~\ref{src:7}, every vertex of $B$ has a core neighbor. A vertex with neighbors in both parity classes would make $D[A]$ semicomplete by~Lemma~\ref{src:15}, so each outside vertex has neighbors in exactly one class. Thus the sets $B_p=\{z\in B:N(z)\cap A_p\ne\varnothing\}$, for $p\in\{0,1\}$, partition $B$. If two adjacent vertices of $B$ each had a neighbor in the same parity class, then $D[A]$ would be semicomplete by Lemma~\ref{src:19}. Both vertices of any pair in $B_p$ have neighbors in $A_p$, so $B_p$ is independent. Moreover, there are no arcs between $B_p$ and $A_{1-p}$, by the definition of $B_p$ and the absence of vertices with neighbors in both core classes. Consequently both $A_0\cup B_1$ and $A_1\cup B_0$ are independent sets. They partition $V(D)$, making $D$ bipartite, a contradiction.

We next prove that $D[B]$ is semicomplete. Since $D[A]$ is semicomplete, we may apply Lemma~\ref{src:13} to any $z\in B$. Either $z$ is adjacent to every vertex of $A$ or it satisfies $\{v_s,\ldots,v_{k+2}\}\strictto z\strictto\{v_0,\ldots,v_r\}$ for some $3\le r<s-1\le k-2$. In either case $z$ has a neighbor in both $A_0$ and $A_1$; in the second case, it dominates $v_0$ and $v_1$. By Lemma~\ref{src:16}, any two distinct vertices of $B$ are adjacent, so $D[B]$ is semicomplete.

Finally, if $z\in B$ has a nonneighbor in $A$, only the second alternative of the classification above is possible. Its final and initial segments contain $v_{k+2}$ and $v_0$, respectively, giving $v_{k+2}\to z\to v_0$.
\end{proof}

In the bipartite case, the missing core adjacencies are confined to an interval within one parity class.

\begin{proposition}\label{prop:bip}
Suppose that $k\ge7$ is odd and $D$ is strong, bipartite, and $k$-quasi-transitive. Choose its bipartition $(X,Y)$ so that $v_i\in X$ exactly when $i$ is even. Then the following hold.
\begin{enumerate}[label=\textup{(\roman*)}]
\item\label{item:complete} The digraphs $D[A]$ and $D[B]$ are semicomplete bipartite with their inherited bipartitions. In particular, the core contains the arcs
\begin{equation}\label{eq:corearcs}
v_i\to v_{i+1}\;(0\le i\le k+1)\;\text{and}\;v_j\to v_i\;(j-i\ge3\text{ odd})\,.
\end{equation}
\item\label{item:boundary} Each $z\in B$ is adjacent to the first two and the last two core vertices in the part opposite to $z$.
\item\label{item:interval} If $z\in B$ is nonadjacent to an opposite-part core vertex $v_i$, its core outneighbors form an initial segment of the opposite-part core vertices, and its core inneighbors form a final segment. These segments are separated by a nonempty interval of nonneighbors, and each contains at least two vertices. In particular, all core outneighbors of $z$ precede $v_i$, all core inneighbors follow $v_i$, and
\begin{equation}\label{eq:bounds}
\begin{cases}
5\le i\le k-2, & z\in X,\quad i\text{ odd},\\
4\le i\le k-3, & z\in Y,\quad i\text{ even}.
\end{cases}
\end{equation}
\end{enumerate}
\end{proposition}
\begin{proof}
For~\ref{item:complete}, the bipartition along $P$ is $A\cap X=A_0$ and $A\cap Y=A_1$. By~Lemma~\ref{src:4}, every pair with one vertex in each set is adjacent, so $D[A]$ is semicomplete bipartite. The forward consecutive arcs belong to $P$. If $j-i\ge3$ is odd, $v_i$ and $v_j$ are adjacent, whereas $v_i\to v_j$ would shorten the $v_i$--$v_j$ subpath of $P$. Thus $v_j\to v_i$, proving~\eqref{eq:corearcs}.

Now take $x\in B\cap X$ and $y\in B\cap Y$. By Lemma~\ref{src:7}, both have core neighbors. Bipartiteness forces $N(x)\cap A_1\ne\varnothing$ and $N(y)\cap A_0\ne\varnothing$. Applying Lemma~\ref{src:16} to the ordered pair $(y,x)$, we conclude that $x$ and $y$ are adjacent. Since this holds for every such pair, $D[B]$ is semicomplete bipartite. If an inherited part is empty, the assertion is vacuous.

For~\ref{item:boundary}, let $z\in B$ and let $A_p$ be its opposite-part core class. By Lemma~\ref{src:7}, $z$ has a neighbor $v_t\in A_p$. If $z\to v_t$ and $v_t$ is not among the first two vertices of $A_p$, then $z\to v_p$ by Lemma~\ref{src:8}. Similarly, if $v_t\to z$ and $v_t$ is not among the last two vertices of $A_p$, then the last vertex of $A_p$ dominates $z$ by the same lemma. Thus $z$ is adjacent to at least one of these four boundary vertices and hence, by Lemma~\ref{src:9}, to all four.

For~\ref{item:interval}, fix $z\in B$ and let $A_p$ be its opposite-part core class. As shown above, $z$ has a neighbor in $A_p$. By Lemma~\ref{src:10}, either $z$ is adjacent to every vertex of $A_p$, or its outneighbors in $A_p$ are exactly the first $r$ vertices and its inneighbors in $A_p$ exactly the last $s$ vertices, where $r,s\ge2$ and $r+s<|A_p|$.

If $z\nonadj v_i$ with $v_i\in A_p$, only the second possibility remains. The two neighbor segments are disjoint and leave a nonempty interval of nonneighbors between them. Since $v_i$ is in this interval, all core outneighbors precede $v_i$ and all core inneighbors follow it. Bipartiteness excludes core neighbors in the other parity class.

If $z\in X$, then $A_p=A_1$ has indices $1,3,\ldots,k,k+2$. Its first two and last two vertices are neighbors, so a nonneighbor has odd index $5\le i\le k-2$. If $z\in Y$, then $A_p=A_0$ has indices $0,2,\ldots,k-1,k+1$, and a nonneighbor has even index $4\le i\le k-3$. These are precisely the bounds in~\eqref{eq:bounds}.
\end{proof}

In either proposition, the core contains the cycle $v_0v_1\cdots v_{k+2}v_0$.

\section{The non-bipartite case}\label{sec:non-bipartite}
Semicompleteness puts the outside vertices having a core nonneighbor on one path. Proposition~\ref{prop:semi} supplies the arcs attaching this path to the core cycle. Every remaining vertex is adjacent to all other vertices, so we can use the following lemma to obtain a Hamiltonian cycle.

\begin{lemma}\label{lem:extension}
Let $D$ be strong, and let $C_0$ be a directed cycle in $D$. If every vertex outside $C_0$ is adjacent to every other vertex of $D$, then $D$ is Hamiltonian.
\end{lemma}
\begin{proof}
Choose a cycle $C$ of maximum order subject to $V(C_0)\subseteq V(C)$, and suppose that $R=V(D)\setminus V(C)$ is nonempty. If $z\in R$ has both an inneighbor and an outneighbor on $C$, it can be inserted between consecutive vertices of $C$. Indeed, the set of vertices of $C$ dominated by $z$ is nonempty. If it is proper, choose a transition along $C$ from its complement into this set. Adjacency supplies the arc from the preceding vertex to $z$. If it is all of $V(C)$, any inneighbor of $z$ and its successor give an insertion. Either case contradicts maximality.

Put $R^+=\{z\in R:V(C)\strictto z\}$ and $R^-=\{z\in R:z\strictto V(C)\}$. Every vertex of $R$ consequently belongs to exactly one of these sets. Both sets are nonempty: if $R^-=\varnothing$, no arc leaves $R$, while if $R^+=\varnothing$, no arc leaves $V(C)$, contradicting strong connectivity. For the same reason, some arc $a\to b$ has $a\in R^+$ and $b\in R^-$; otherwise no arc leaves $R^+$. Replacing any arc $c\to c^+$ of $C$ by $c\to a\to b\to c^+$ gives a larger cycle, a contradiction.
\end{proof}

\begin{proof}[Proof of Theorem~\ref{thm:odd}]
Fix~\eqref{eq:geodesic}, and let $M=\{z\in B:z\nonadj v_i\text{ for some }v_i\in A\}$. By Proposition~\ref{prop:semi}, $D[M]$ is semicomplete. Every nonempty semicomplete digraph has a Hamiltonian path, by induction on its order. Given a path through all but one vertex $z$, insert $z$ before the first vertex it dominates, or append it if there is no such vertex. If the insertion point has a predecessor, that predecessor dominates $z$ by adjacency and the choice of the insertion point.

If $M\ne\varnothing$, choose a path $z_1\cdots z_m$ spanning $M$. The cycle $v_0v_1\cdots v_{k+2}z_1\cdots z_mv_0$ then exists by Proposition~\ref{prop:semi}. If $M=\varnothing$, use the core cycle. In either case there is a cycle $C_0$ with vertex set $A\cup M$. Every vertex of $B\setminus M$ is adjacent to all vertices of $A$ by definition and to all other vertices of $B$ by Proposition~\ref{prop:semi}. Applying Lemma~\ref{lem:extension} to $C_0$, we conclude that $D$ is Hamiltonian.
\end{proof}

\section{Separation in the bipartite case}\label{sec:separation}
Throughout this section, $k\ge7$ is odd and $D$ is strong, bipartite, and $k$-quasi-transitive, with $\diam(D)\ge k+2$. Fix~\eqref{eq:geodesic} and the notation of Proposition~\ref{prop:bip}. Our aim is to show that every $v_0$--$z$ path and every $z$--$v_{k+2}$ path passes through $v_i$ whenever $z\in B$ and $v_i\in A$ are opposite-part nonneighbors. The proofs combine outside path segments with core paths of prescribed length to contradict $k$-quasi-transitivity.

The arcs in~\eqref{eq:corearcs} allow forward motion along consecutive indices and backward motion across any odd index difference of at least three. In the next lemma we use these arcs to construct simple paths of prescribed length.

\begin{lemma}\label{lem:paths}
Put $n=k+3$. If $0\le t<s\le n-1$, $3\le L\le n-4$, and $L\equiv s-t\pmod2$, then the arcs in~\eqref{eq:corearcs} contain a simple $v_s$--$v_t$ path of length $L$.
\end{lemma}
\begin{proof}
Set $m=n-L$ and $r=s-t$. Then $m\ge4$ and $m\equiv r\pmod2$, since $n$ is even. All intervals of indices below are traversed in increasing order.

Suppose first that $r\le m$. Put $b=\max\{0,t-L+1\}$ and $a=b+L+r-1$. Consider
\begin{equation}\label{eq:path-two}
v_sv_{s+1}\cdots v_a\,v_bv_{b+1}\cdots v_t\,.
\end{equation}
The choice of $b$ gives $0\le b\le t$ and $b+L-1\ge t$, whence $s\le a$. Also $b\le m-r$: this follows from $r\le m$ when $b=0$, and from $t+r=s\le n-1$ when $b=t-L+1$. Hence $a\le n-1$. Thus both intervals lie within the core and are disjoint, since $t<s$. Their joining arc exists because $a-b=L+r-1$ is odd and at least three. The length of~\eqref{eq:path-two} is $(a-s)+1+(t-b)=L$.

Now suppose that $r>m$. The parity relation gives $r\ge m+2$. Let $u$ be the smallest even integer strictly greater than $t+1$, and put $q=u+r-m-1$. Consider
\begin{equation}\label{eq:path-three}
v_sv_{s+1}\cdots v_{n-1}\,v_uv_{u+1}\cdots v_q\,v_0v_1\cdots v_t\,.
\end{equation}
The middle interval contains $r-m\ge2$ vertices. Since $u$ is even and $r-m$ is even, $q$ is odd. Also, $t+2\le u\le t+3\le s-m+1\le n-m$ and $3\le q\le s-m+2\le s-2$. Thus the three intervals are pairwise disjoint. The backward differences $(n-1)-u$ and $q-0$ are odd and at least three, so both joining arcs exist. The number of vertices in~\eqref{eq:path-three} is $(n-s)+(r-m)+(t+1)=L+1$, as required.
\end{proof}

Fix an outside vertex $z$ and an opposite-part core nonneighbor $v_i$. The core paths just constructed restrict the arcs incident with $z$ and prevent an outside path from entering the core beyond $v_i$.

\begin{lemma}\label{lem:local}
Let $z\in B$ and $v_i\in A$ lie in opposite parts, with $z\nonadj v_i$, and let $w\in B$ lie in the part opposite to $z$.
\begin{enumerate}[label=\textup{(\roman*)}]
\item If $t<i$, $v_t$ is in the same part as $z$, and $v_t\to w$, then $w\not\to z$ and $z\to w$.
\item If $t>i$, $v_t$ is in the same part as $z$, and $w\to v_t$, then $z\not\to w$ and $w\to z$.
\end{enumerate}
\end{lemma}
\begin{proof}
For the first assertion, suppose that $w\to z$. Since $k-2\in[3,k-1]$ and both $k-2$ and $i-t$ are odd, we may apply Lemma~\ref{lem:paths} to obtain a core $v_i$--$v_t$ path of length $k-2$. Appending $w,z$ gives a simple path of length $k$ with nonadjacent endpoints $v_i,z$, a contradiction. Hence $w\not\to z$, and the semicomplete bipartite structure of $D[B]$ gives $z\to w$.

For the second assertion, suppose that $z\to w$. Since $t>i$ and $t-i$ is odd, we obtain a core $v_t$--$v_i$ path of length $k-2$ by Lemma~\ref{lem:paths}. Prepending $z,w$ gives a simple path of length $k$ from $z$ to its nonneighbor $v_i$, a contradiction. Therefore $z\not\to w$, and opposite-part adjacency in $D[B]$ gives $w\to z$.
\end{proof}

\begin{lemma}\label{lem:outside}
Let $z\in B$ and $v_i\in A$ lie in opposite parts, with $z\nonadj v_i$. There is no $z$--$v_t$ path with $t>i$ whose internal vertices all lie in $B$.
\end{lemma}
\begin{proof}
Suppose that such a path exists, and choose one of minimum length. Write it as $Q=w_0w_1\cdots w_rv_t$, where $w_0=z$, all $w_j$ lie in $B$, and $t>i$. Let $\ell=r+1$ be its length. Proposition~\ref{prop:bip}\ref{item:interval} excludes $z\to v_t$ when $v_t$ is opposite to $z$, and bipartiteness excludes it otherwise. Thus $\ell\ge2$.

Set $h=2$ if $z\in X$ and $h=1$ if $z\in Y$, so that $v_h$ is an early core vertex in the same part as $z$. Define $j=k-1-i+h$, so that $j+1+i-h=k$. This choice will allow us to complete a sufficiently long outside path to a path of length $k$. If $z\in X$,~\eqref{eq:bounds} gives odd $i\in[5,k-2]$ and $j=k+1-i$; if $z\in Y$, it gives even $i\in[4,k-3]$ and $j=k-i$. In either case, $j$ is odd, $3\le j\le k-4$, and $h<i$.

Suppose first that $r\ge j$. The vertex $w_j$ is opposite to $z$ and to $v_h$. Minimality of $Q$ excludes $z\to w_j$, since this arc would shorten $Q$. Thus $w_j\to z$ by Proposition~\ref{prop:bip}\ref{item:complete}. If $v_h\to w_j$, then $w_j\not\to z$ by Lemma~\ref{lem:local}(i), a contradiction. Since $v_h$ is one of the first two core vertices in the part opposite to $w_j$, Proposition~\ref{prop:bip}\ref{item:boundary} ensures that $w_j$ and $v_h$ are adjacent. Hence $w_j\to v_h$. The path $z,w_1,\ldots,w_j,v_h,v_{h+1},\ldots,v_i$ is simple and has length $j+1+i-h=k$, contradicting $z\nonadj v_i$.

Consequently $r<j$, so $2\le\ell=r+1\le j\le k-4$. Put $L=k-\ell$. Then $4\le L\le k-2$. Since $z$ and $v_i$ are opposite, bipartiteness of $Q$ implies $t-i\equiv1-\ell\equiv k-\ell\pmod2$. We can therefore apply Lemma~\ref{lem:paths} to obtain a core $v_t$--$v_i$ path of length $L$. Its only intersection with $Q$ is $v_t$, so their concatenation is a simple path of length $k$ from $z$ to $v_i$, again a contradiction.
\end{proof}

\begin{theorem}\label{thm:separator}
Let $z\in B$ and $v_i\in A$ lie in opposite parts, with $z\nonadj v_i$. Then $D-v_i$ contains neither a $v_0$--$z$ path nor a $z$--$v_{k+2}$ path.
\end{theorem}
\begin{proof}
Suppose that there is a $z$--$v_{k+2}$ path $Q$ avoiding $v_i$. The initial segment of $Q$ up to its first core vertex has all internal vertices in $B$. Since $Q$ avoids $v_i$, it follows from Lemma~\ref{lem:outside} that this core vertex has index less than $i$. Let $v_b$ be the first core vertex on $Q$ with $b>i$, and let $v_a$ be the preceding core vertex on $Q$. Thus $a<i<b$, and the intervening subpath is
\begin{equation}\label{eq:crossing}
v_a,w_1,\ldots,w_r,v_b, \text{where } w_1,\ldots,w_r\in B\,.
\end{equation}
The shortest-path property of $P$ excludes the direct arc $v_a\to v_b$, so $r\ge1$. None of the $w_j$ equals $z$, since $Q$ is simple and starts at $z$.

If $v_a$ lies in the same part as $z$, we may apply Lemma~\ref{lem:local}(i) to the arc $v_a\to w_1$ to obtain $z\to w_1$. The path $z,w_1,\ldots,w_r,v_b$ contradicts Lemma~\ref{lem:outside}.

Suppose instead that $v_a$ lies in the same part as $v_i$. If $r=1$, then $v_a$ and $v_b$ are in the same part, so $a,i,b$ have the same parity. Hence $b-a\ge4$, whereas~\eqref{eq:crossing} has length two, contradicting $d(v_a,v_b)=b-a$. Thus $w_2$ exists and is opposite to $z$.

If $w_2\to z$, use Lemma~\ref{lem:paths} to obtain a core $v_i$--$v_a$ path of length $k-3$. Indeed, $4\le k-3\le k-1$, and $i-a$ and $k-3$ are both even. Appending $w_1,w_2,z$ gives a simple path of length $k$ with nonadjacent endpoints. Therefore $w_2\not\to z$, and opposite-part adjacency in $D[B]$ implies $z\to w_2$. Now $z,w_2,\ldots,w_r,v_b$ contradicts Lemma~\ref{lem:outside}.

This proves that $D-v_i$ has no $z$--$v_{k+2}$ path. Reverse all arcs of $D$ and use the reversed shortest path $v_{k+2}\cdots v_0$. Strong connectivity, bipartiteness, $k$-quasi-transitivity, and nonadjacency are preserved, and the reversed path is again shortest of length $k+2$. The assertion just proved excludes a $z$--$v_0$ path in the converse digraph with $v_i$ deleted, equivalently a $v_0$--$z$ path in $D-v_i$.
\end{proof}

\begin{corollary}\label{cor:cycle}
If a directed cycle contains $z\in B$ and at least one vertex of $A$, then it contains every opposite-part core vertex nonadjacent to $z$.
\end{corollary}
\begin{proof}
Suppose that the cycle avoids such a vertex $v_i$, and choose a core vertex $v_a$ on the cycle. If $a<i$, follow $P$ from $v_0$ to $v_a$ and then the cycle to $z$. This is a $v_0$--$z$ walk avoiding $v_i$, and hence contains a path of that kind. If $a>i$, follow the cycle from $z$ to $v_a$ and then $P$ to $v_{k+2}$. This yields a $z$--$v_{k+2}$ walk avoiding $v_i$. Both cases contradict Theorem~\ref{thm:separator}.
\end{proof}

\section{Merging cycles}\label{sec:merging}
The merging arguments below require complete opposite-part adjacency between distinct cycles, but impose no such condition within a cycle. We first join two cycles by a two-arc exchange whenever one is available. If no such exchange exists, the directions of the cross arcs force a four-arc exchange.

\begin{figure}[tb]
\centering
\begin{tikzpicture}[x=1cm,y=1cm,font=\small]
\node[font=\small] at (7,3.05) {(a) Cutting the cycles};
\node[cv,label=above:$c_1$] (u1) at (1,5.7) {};
\node[cv,label=above:$c_2$] (u2) at (5,5.7) {};
\node[cv,label=below:$c_3$] (u3) at (5,4.1) {};
\node[cv,label=below:$c_0$] (u0) at (1,4.1) {};
\draw[keepC] (u1)--node[above] {$c_1Cc_2$}(u2);
\draw[removed] (u2)--(u3);
\draw[keepC] (u3) to[bend left=13] node[below] {$c_3Cc_0$}(u0);
\draw[removed] (u0)--(u1);
\node[text=cycleC] at (3,4.95) {cycle $C$};
\node[ev,label=above:$e_{a+1}$] (w1) at (9,5.7) {};
\node[ev,label=above:$e_b$] (w2) at (13,5.7) {};
\node[ev,label=below:$e_{b+1}$] (w3) at (13,4.1) {};
\node[ev,label=below:$e_a$] (w0) at (9,4.1) {};
\draw[keepE] (w1) to[bend left=8] node[above] {$e_{a+1}Ee_b$}(w2);
\draw[removed] (w2)--(w3);
\draw[keepE] (w3) to[bend left=13] node[below] {$e_{b+1}Ee_a$}(w0);
\draw[removed] (w0)--(w1);
\node[text=cycleE] at (11,4.95) {cycle $E$};
\node[font=\small] at (7,-1.25) {(b) Reconnecting the paths};
\node[cv,label=above:$c_1$] (c1) at (1,1.25) {};
\node[cv,label=above:$c_2$] (c2) at (4,1.25) {};
\node[ev,label=above:$e_{a+1}$] (ea1) at (8,1.25) {};
\node[ev,label=above:$e_b$] (eb) at (12.5,1.25) {};
\node[cv,label=below:$c_3$] (c3) at (12.5,-.55) {};
\node[cv,label=below:$c_0$] (c0) at (9,-.55) {};
\node[ev,label=below:$e_{b+1}$] (eb1) at (5,-.55) {};
\node[ev,label=below:$e_a$] (ea) at (1,-.55) {};
\draw[keepC] (c1)--node[above] {$C$}(c2);
\draw[added] (c2)--node[badge] {1}(ea1);
\draw[keepE] (ea1)--node[above] {$E$}(eb);
\draw[added] (eb)--node[badge] {2}(c3);
\draw[keepC] (c3)--node[below] {$C$}(c0);
\draw[added] (c0)--node[badge] {3}(eb1);
\draw[keepE] (eb1)--node[below] {$E$}(ea);
\draw[added] (ea)--node[badge] {4}(c1);
\node[font=\footnotesize,text=black!65] at (6.8,.35) {Every original vertex is retained; only the four connecting arcs change.};
\draw[removed] (.5,-1.9)--(1.3,-1.9); \node[anchor=west,font=\footnotesize] at (1.4,-1.9) {deleted arc};
\draw[added] (4.6,-1.9)--(5.4,-1.9); \node[anchor=west,font=\footnotesize] at (5.5,-1.9) {new cross arc};
\draw[keepC] (9.3,-1.9)--(10.1,-1.9); \node[anchor=west,font=\footnotesize] at (10.2,-1.9) {retained path};
\end{tikzpicture}
\caption{Four-arc exchange. Solid paths are retained; dash-dotted arcs are inserted.}
\label{fig:exchange}
\end{figure}
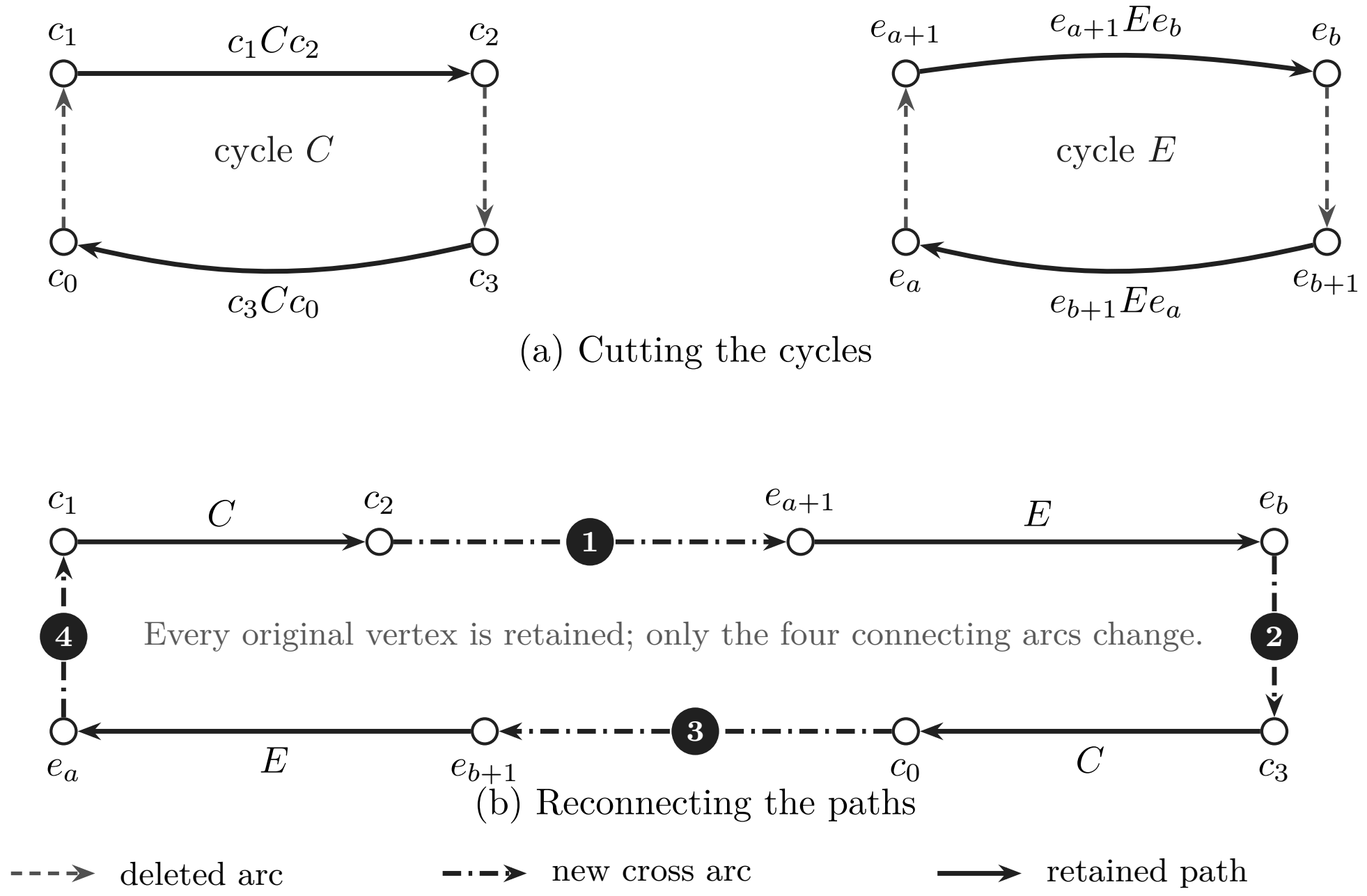

\begin{lemma}\label{lem:two-cycles}
Let $C,E$ be vertex-disjoint directed cycles in a bipartite digraph. Suppose that every opposite-part pair with one vertex on each cycle is adjacent. If there is an arc from $C$ to $E$ and an arc from $E$ to $C$, then the subdigraph induced by $V(C)\cup V(E)$ is Hamiltonian.
\end{lemma}
\begin{proof}
Write $C=c_0c_1\cdots c_{2p-1}c_0$ and $E=e_0e_1\cdots e_{2q-1}e_0$, choosing the indices so that same-parity indices on the two cycles represent vertices in the same part. Subscripts on each cycle are read modulo its length.

If $c_i\to e_j$ and $e_{j-1}\to c_{i+1}$, replace the arcs $c_i\to c_{i+1}$ and $e_{j-1}\to e_j$ by these two cross arcs. The resulting cycle is $c_{i+1}Cc_i e_jEe_{j-1}c_{i+1}$. We may therefore assume that no such pair exists.

Under this assumption, an arc $c_i\to e_j$ excludes $e_{j-1}\to c_{i+1}$. Opposite-part adjacency then forces $c_{i+1}\to e_{j-1}$. Iterate the permutation $T(i,j)=(i+1,j-1)$ on opposite-parity index pairs. Since the orbit is finite, iteration returns to the starting pair. Thus, if one pair in the orbit supports an arc from $C$ to $E$, every pair does. At each step the forbidden two-arc exchange excludes the reverse arc at the next pair, so no reverse arc occurs anywhere on the orbit. If there is no arc from $C$ to $E$ on an orbit, adjacency gives the strict direction from $E$ to $C$ throughout it.

Put $g=\gcd(2p,2q)$. Two pairs $(i,j)$ and $(i',j')$ lie on the same orbit exactly when $i+j\equiv i'+j'\pmod g$. Indeed, reaching the second pair after $m$ steps requires $m\equiv i'-i\pmod{2p}$ and $m\equiv j-j'\pmod{2q}$, and these congruences are compatible exactly under the stated condition. Only odd residues occur. For an odd residue $s$ modulo $g$, let $f(s)=1$ if its orbit is directed from $C$ to $E$, and let $f(s)=0$ otherwise. Arguments of $f$ are read modulo $g$.

Both directions occur by hypothesis, so the cyclic sequence $f(1),f(3),\ldots,f(g-1)$ is nonconstant. In particular, $g\ge4$ and both cycles have length at least four. Choose a transition from $0$ to $1$ and a transition from $1$ to $0$. They give distinct even representatives $a,b\in\{0,2,\ldots,g-2\}$ satisfying $f(a+1)=0$, $f(a+3)=1$, $f(b+1)=1$, and $f(b+3)=0$. Since $g\le2q$, these representatives are also indices on $E$. The index sums of the following four pairs are $a+3$, $b+3$, $b+1$, and $a+1$, respectively. The chosen values of $f$ therefore give the arcs $c_2\to e_{a+1}$, $e_b\to c_3$, $c_0\to e_{b+1}$, and $e_a\to c_1$.

Consequently,
\begin{equation}\label{eq:four-exchange}
c_1Cc_2\,e_{a+1}Ee_b\,c_3Cc_0\,e_{b+1}Ee_a\,c_1
\end{equation}
is a directed cycle. The two $C$ segments partition $V(C)$. Deleting the distinct arcs $e_a\to e_{a+1}$ and $e_b\to e_{b+1}$ partitions $E$ into the two displayed $E$ segments. Thus~\eqref{eq:four-exchange} visits every vertex of the union exactly once; see Figure~\ref{fig:exchange}.
\end{proof}

\begin{lemma}\label{lem:factor-merge}
Let $H$ be a strong bipartite digraph with a cycle factor $\mathcal F$. Suppose that every opposite-part pair belonging to distinct cycles of $\mathcal F$ is adjacent. Then $H$ is Hamiltonian.
\end{lemma}
\begin{proof}
Let $(X,Y)$ be a bipartition of $H$. Among the factors obtainable by merging cycles of $\mathcal F$, choose one with the fewest cycles. The assumed cross adjacency is preserved by merging, since vertices on distinct resulting cycles belonged to distinct original cycles. If two remaining cycles had arcs in both directions between them, we could merge them by Lemma~\ref{lem:two-cycles}, contrary to the choice of the factor. Hence the cross arcs between any two remaining cycles are all directed the same way. Each cycle contains vertices in both parts, so at least one cross arc exists between every pair of cycles. Contracting the cycles therefore gives a tournament. This tournament is strong because every directed path in $H$ projects to a directed walk between the contracted vertices.

If more than one cycle remains, the strong quotient tournament has a directed cycle. List its corresponding factor cycles in cyclic order as $D_1,\ldots,D_m$. On each $D_j$, choose an arc $x_j\to y_j$ with $x_j\in X$ and $y_j\in Y$. The quotient direction and complete opposite-part adjacency supply every arc $x_j\to y_{j+1}$, with subscripts modulo $m$. Delete $x_j\to y_j$ from each selected cycle and add these cross arcs. Deleting $x_j\to y_j$ leaves the path $y_jD_jx_j$ through all vertices of $D_j$. The new arcs join these disjoint paths in the cyclic order $D_1,\ldots,D_m$, producing one directed cycle. This contradicts the minimal number of cycles.
\end{proof}

\section{The cycle-factor characterization}\label{sec:factor}
\begin{proof}[Proof of Theorem~\ref{thm:factor}]
A Hamiltonian cycle is a cycle factor, so it remains to prove sufficiency. Suppose that $D$ has a cycle factor. If $k$ is even, then $k\ge6$, and~\cite[Theorem 15]{HJ} gives pancyclicity and hence Hamiltonicity. If $k$ is odd and $D$ is non-bipartite, apply Theorem~\ref{thm:odd}.

It remains to consider odd $k\ge7$ and bipartite $D$. Use the notation of Section~\ref{sec:separation}. Let $M$ be the set of vertices in $B$ having an opposite-part nonneighbor in $A$, and put $U=B\setminus M$.

Every vertex of $U$ is adjacent to every opposite-part vertex of $D$: adjacency to $A$ follows from its definition, and adjacency to $B$ follows from Proposition~\ref{prop:bip}. We first concentrate the core on one factor cycle, and then absorb all remaining cycles meeting $M$.

\begin{claim}\label{lem:concentrate}
If $D$ has a cycle factor, it has a cycle factor with one cycle containing all of $A$.
\end{claim}
\begin{proof}[Proof of Claim \ref{lem:concentrate}]\renewcommand{\qedsymbol}{$\square$}
Take a cycle factor $\mathcal F$, and let $\mathcal F_A$ consist of its cycles meeting $A$. Every opposite-part pair in distinct cycles of $\mathcal F_A$ is adjacent. Opposite-part pairs within $A$ or within $B$ are adjacent by Proposition~\ref{prop:bip}. The only other possibility is a pair $v_i\in A$, $z\in B$. If they were nonadjacent, the cycle containing $z$ meets $A$, and Corollary~\ref{cor:cycle} would force it to contain $v_i$ as well, contrary to their belonging to distinct factor cycles.

Let $W$ be the union of the vertex sets of the cycles in $\mathcal F_A$. The cycle $v_0v_1\cdots v_{k+2}v_0$ lies in $D[A]$, so $D[A]$ is strong. Every vertex of $W$ lies on a cycle meeting $A$, so it can reach and be reached from $A$ within $D[W]$. Thus $D[W]$ is strong. Applying Lemma~\ref{lem:factor-merge} to $D[W]$ and its cycle factor $\mathcal F_A$, we obtain a Hamiltonian cycle of $D[W]$ containing all of $A$. Together with the unchanged cycles of $\mathcal F\setminus\mathcal F_A$, this cycle gives the required factor.
\end{proof}

The core vertices may occur on the enlarged cycle in a different order from their order on $P$. The separation theorem identifies a segment between $v_{k+2}$ and $v_0$ on which an outside cycle can be inserted.

\begin{claim}\label{lem:absorb}
Let $C,E$ be vertex-disjoint directed cycles such that $A\subseteq V(C)$, $V(E)\subseteq B$, and $V(E)\cap M\ne\varnothing$. Then $D[V(C)\cup V(E)]$ is Hamiltonian.
\end{claim}
\begin{proof}[Proof of Claim \ref{lem:absorb}]\renewcommand{\qedsymbol}{$\square$}
Write $v^*=v_{k+2}$. We first prove
\begin{equation}\label{eq:strict-boundary}
v^*\strictto(V(E)\cap X)\text{ and }(V(E)\cap Y)\strictto v_0\,.
\end{equation}
Choose $z\in V(E)\cap M$ and an opposite-part core nonneighbor $v_i$. Each $e\in V(E)\cap X$ is adjacent to $v^*$ by Proposition~\ref{prop:bip}\ref{item:boundary}. An arc $e\to v^*$, preceded by the segment of $E$ from $z$ to $e$ (of length zero if $e=z$), would give a $z$--$v^*$ path avoiding $v_i$, contrary to Theorem~\ref{thm:separator}. Hence $v^*\to e$ and $e\not\to v^*$. An arc $v_0\to e$ with $e\in V(E)\cap Y$, followed by the segment of $E$ from $e$ to $z$, would give a $v_0$--$z$ path avoiding $v_i$, contradicting the other assertion of Theorem~\ref{thm:separator}. Boundary adjacency supplies $e\to v_0$, proving~\eqref{eq:strict-boundary}.

Let $J$ be the set of core vertices having an opposite-part nonneighbor on $E$. For each $v_i\in J$, choose such a nonneighbor $z_i\in V(E)$. There is no $v_0$--$v^*$ path in $D-v_i$. Otherwise follow such a path by an arc from $v^*$ to any vertex of $V(E)\cap X$, and then follow $E$ to $z_i$. This produces a $v_0$--$z_i$ walk avoiding $v_i$, and hence a path, contrary to Theorem~\ref{thm:separator}. The walk argument applies even if the initial path intersects $E$.

It follows that the segment $v_0Cv^*$ contains all of $J$. Neither endpoint belongs to $J$, by Proposition~\ref{prop:bip}\ref{item:boundary}. The complementary segment $Q=v^*Cv_0=q_0q_1\cdots q_r$ therefore avoids $J$. Every vertex of $Q$ is adjacent to every opposite-part vertex of $E$. For a core vertex this follows from the definition of $J$, and for a vertex in $B$ it follows from the semicomplete bipartite structure of $D[B]$.

If an arc $q_j\to q_{j+1}$ of $Q$ and an arc $e\to e^+$ of $E$ satisfy
\begin{equation}\label{eq:insert}
q_j\to e^+\text{ and } e\to q_{j+1}\,,
\end{equation}
replace the two original arcs by these cross arcs. This merges $C$ and $E$, since the path $e^+Ee$ contains all vertices of $E$; see Figure~\ref{fig:absorb}. In the figure, solid arrows represent retained segments of $C$ and $E$. Either end segment of $Q$ may have length zero, in which case its endpoints coincide.

Suppose no such pair exists. We propagate domination along $Q$ to obtain a contradiction at $v_0$. By~\eqref{eq:strict-boundary}, $q_0=v^*$ strictly dominates every opposite-part vertex of $E$. Suppose inductively that $q_j$ has this property, where $j<r$. For any $e\in V(E)$ opposite to $q_{j+1}$, its successor $e^+$ is opposite to $q_j$, so $q_j\to e^+$. The failure of~\eqref{eq:insert} implies $e\not\to q_{j+1}$, and cross adjacency implies $q_{j+1}\to e$. Induction reaches $q_r=v_0$ and excludes every arc from $V(E)\cap Y$ to $v_0$, contradicting~\eqref{eq:strict-boundary}. An insertion pair must therefore exist.
\end{proof}

\begin{figure}[tb]
\centering
\begin{tikzpicture}[x=1cm,y=1cm,font=\small]
\node[font=\small] at (3.35,1.65) {(a) A segment avoiding $J$};
\node[cv,label=left:$v_0$] (a0) at (.9,4.1) {};
\node[cv,label=right:$v^*$] (as) at (5.6,4.1) {};
\draw[keepC] (a0) .. controls (1.6,6) and (4.9,6) .. node[above,font=\footnotesize] {$v_0Cv^*$ contains all of $J$} (as);
\draw[keepC] (as) .. controls (4.9,2.2) and (1.6,2.2) .. node[below,font=\footnotesize] {$Q=v^*Cv_0$ avoids $J$} (a0);
\node[text=cycleC,font=\large] at (3.25,4.1) {$C$};
\node[font=\small] at (10.5,1.65) {(b) Exchanging two arcs};
\node[ev,label=left:$e^+$] (be1) at (8.5,4.85) {};
\node[ev,label=right:$e$] (be0) at (12.5,4.85) {};
\draw[keepE] (be1) to[bend left=20] node[above,font=\footnotesize] {$e^+Ee$}(be0);
\draw[removed] (be0) to[bend left=14] (be1);
\node[cv,label=below:$q_j$] (bq0) at (8.5,2.9) {};
\node[cv,label=below:$q_{j+1}$] (bq1) at (12.5,2.9) {};
\draw[removed] (bq0)--node[below,font=\footnotesize] {an arc of $Q$}(bq1);
\draw[added] (bq0)--(be1);
\draw[added] (be0)--(bq1);
\node[font=\small] at (7,-2.95) {(c) The merged cycle};
\node[cv,label=left:$v^*$] (vs) at (.9,-.65) {};
\node[cv,label=below:$q_j$] (qj) at (4.2,-.65) {};
\node[cv,label=below:$q_{j+1}$] (qn) at (9.8,-.65) {};
\node[cv,label=right:$v_0$] (v0) at (13.1,-.65) {};
\node[ev,label=above:$e^+$] (ep) at (4.2,.25) {};
\node[ev,label=above:$e$] (e) at (9.8,.25) {};
\draw[keepC] (vs)--node[above] {$Q$}(qj);
\draw[added] (qj)--(ep);
\draw[keepE] (ep)--node[above,font=\footnotesize] {all vertices of $E$}(e);
\draw[added] (e)--(qn);
\draw[keepC] (qn)--node[above] {$Q$}(v0);
\draw[keepC] (v0) .. controls (11.2,-2.35) and (2.8,-2.35) .. node[below,font=\footnotesize] {$v_0Cv^*$ is retained unchanged}(vs);
\end{tikzpicture}
\caption{Cycle absorption. Dashed arcs are deleted; dash-dotted arcs are inserted.}
\label{fig:absorb}
\end{figure}
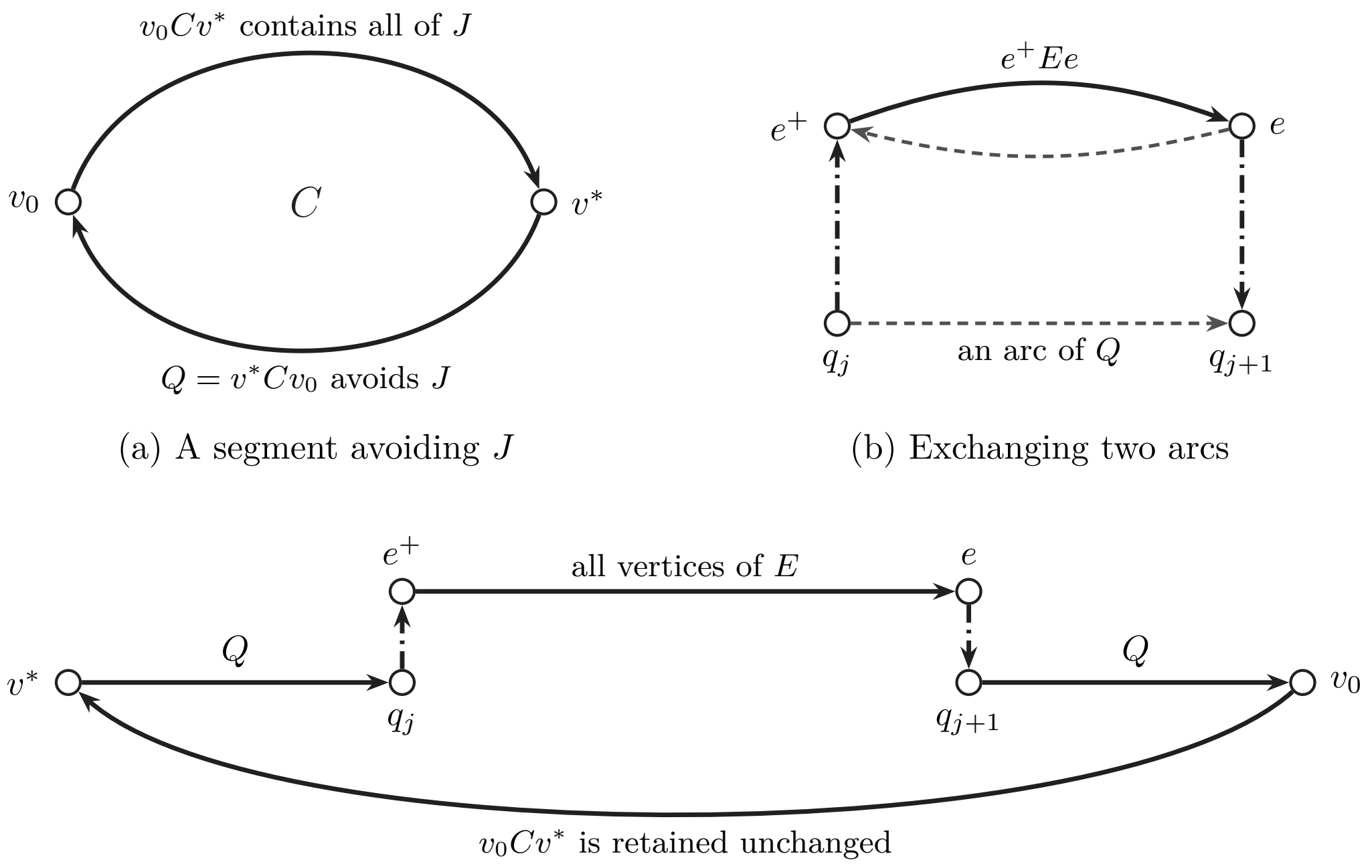

Claim~\ref{lem:concentrate} gives a factor with one cycle $C$ containing all of $A$. Every other factor cycle lies in $B$. Whenever such a cycle meets $M$, merge it into $C$ using Claim~\ref{lem:absorb}. Each merge retains all vertices already on $C$, so the process terminates with one factor cycle containing $A\cup M$ and all other factor cycles contained in $U$.

Every opposite-part pair on distinct cycles of the resulting factor is adjacent: at least one of its vertices belongs to $U$, and such a vertex is adjacent to every opposite-part vertex of $D$. Since $D$ is strong, it has a Hamiltonian cycle by Lemma~\ref{lem:factor-merge}.
\end{proof}

\section*{Declaration of AI use}
AI tools were used to assist in developing the proofs of Propositions \ref{prop:semi} and \ref{prop:bip}, as well as
with literature searches, mathematical discussions, manuscript drafting and revision, and LaTeX
preparation. The authors retain full responsibility for the mathematical correctness, originality,
and final content of the paper.


\begin{thebibliography}{99}
\bibitem{AH} J. Alva-Samos and C. Hern\'andez-Cruz, $k$-quasi-transitive digraphs of large diameter, \emph{Journal of Graph Theory} \textbf{96} (2021), 310--325.

\bibitem{BG} J. Bang-Jensen and G. Gutin, Vertex heaviest paths and cycles in quasi-transitive digraphs, \emph{Discrete Mathematics} \textbf{163} (1997), 217--223.

\bibitem{BGH} J. Bang-Jensen, G. Gutin and J. Huang, A sufficient condition for a semicomplete multipartite digraph to be Hamiltonian, \emph{Discrete Mathematics} \textbf{161} (1996), 1--12.

\bibitem{BH} J. Bang-Jensen and J. Huang, Quasi-transitive digraphs, \emph{Journal of Graph Theory} \textbf{20} (1995), 141--161.

\bibitem{BN} J. Bang-Jensen and M. H. Nielsen, Minimum cycle factors in quasi-transitive digraphs, \emph{Discrete Optimization} \textbf{5} (2008), 121--137.

\bibitem{CGG} P. Cano, H. Galeana-S\'anchez and I. A. Goldfeder, Some results on the existence of Hamiltonian cycles in $\mathcal P$-compositions of bipartite digraphs, \emph{AKCE International Journal of Graphs and Combinatorics} \textbf{17} (2020), 713--719.

\bibitem{GGU} H. Galeana-S\'anchez, I. A. Goldfeder and I. Urrutia, On the structure of strong $3$-quasi-transitive digraphs, \emph{Discrete Mathematics} \textbf{310} (2010), 2495--2498.

\bibitem{HM} R. H\"aggkvist and Y. Manoussakis, Cycles and paths in bipartite tournaments with spanning configurations, \emph{Combinatorica} \textbf{9} (1989), 33--38.

\bibitem{HG} C. Hern\'andez-Cruz and H. Galeana-S\'anchez, $k$-kernels in $k$-transitive and $k$-quasi-transitive digraphs, \emph{Discrete Mathematics} \textbf{312} (2012), 2522--2530.

\bibitem{HJ} C. Hern\'andez-Cruz and M. A. Ju\'arez-Camacho, Pancyclicity in strong $k$-quasi-transitive digraphs of large diameter, \emph{Applied Mathematics and Computation} \textbf{408} (2021), 126319.

\bibitem{WZ} R. Wang and H. Zhang, Hamiltonian paths in $k$-quasi-transitive digraphs, \emph{Discrete Mathematics} \textbf{339} (2016), 2094--2099.

\bibitem{WZodd} R. Wang and H. Zhang, The structure of strong $k$-quasi-transitive digraphs with large diameters,  arXiv: 2006.06333, 2020.
\end{thebibliography}
\end{document}